\documentclass[12pt]{amsart}
\usepackage{amssymb}
\usepackage{graphicx}
\newtheorem{Def}{Definition}

\newtheorem{Thm}{Theorem}
\newtheorem{Cor}{Corollary}
\newtheorem{Rem}{Remark}
\newenvironment{Pf}{ Proof.}{\(\square\)}

\title[A geometric application of Lagrange multipliers...]{A geometric application of Lagrange multipliers: extremal compatible linear connections}
\author{Csaba Vincze}
\address{Inst. of Math., Univ. of Debrecen \\
H-4002 Debrecen, P.O.Box 400 \\
Hungary}
\email{csvincze@science.unideb.hu}
\author{M\'{a}rk Ol\'{a}h}
\address{Institute of Mathematics, University of Debrecen, H-4002 Debrecen, P. O. Box 400, Hungary \newline
\indent HUN-REN-DE Equations, Functions, Curves and their Applications Research Group}
\email{olah.mark@science.unideb.com}
\keywords{Extremal linear connections, Lagrange multipliers on function spaces, Finsler spaces, Generalized Berwald spaces, Intrinsic Geometry}
\subjclass{53C60, 58B20}
\thanks{The research is supported by HUN-REN-DE Equations, Functions, Curves and their Applications Research Project.}
\begin{document}
\begin{abstract}
It is well-known that the L\'{e}vi-Civita connection of a Riemannian manifold is a metric linear connection which is uniquely determined by the vanishing of its torsion. The L\'{e}vi-Civita connection is extremal in the sense that it is a metric linear connection minimizing the length of its torsion point by point. We can also investigate the idea of the extremal linear connection in case of a more general class of spaces with more general (not necessarily quadratic) indicatrix hypersurfaces in the tangent spaces. 

Beyond Riemannian spaces the existence of metric or, in other words, compatible linear connections on the base manifold is far from being automatic. Another problem is the intrinsic characterization of the extremal one including the solution of the existence problem. The first step is to provide the Riemann metrizability of the compatible linear connections. The Riemannian environment establishes a one-to-one correspondence between the linear connections and their torsions, and a relatively simple way of measuring the length of the torsion tensors is also given. The second step is to solve a hybrid conditional extremum problem at each point of the base manifold all of whose constraint equations (compatibility equations) involve functions defined on the indicatrix hypersurface. The objective function to be minimized is a quadratic squared norm function defined on the finite dimensional fiber (vector space) of the torsion tensor bundle. The paper is devoted to the expression of the solution by using the method of Lagrange multipliers on function spaces point by point. We present a necessary and sufficient condition of the solvability and the solution is also given in terms of intrinsic quantities affecting the uniform size of the linear isometry groups of the indicatrices. This completes the description of differential geometric spaces admitting compatible linear connections on the base manifold. They are called generalized Berwald spaces in Finsler geometry.

The presented method of extremals can be applied to the solution of each problem of similar type: find Riemann metrizable linear connections such that the parallel transports preserve a smoothly varying family of quantities in the tangent spaces of the base manifold. 
\end{abstract}

\maketitle

\section*{Introduction}

Finsler metrics are constituted by smoothly varying families of Minkowski functionals (norm functions in the centrally symmetric case) in the tangent spaces of the base manifold. The Minkowski functionals allow us to measure the length of curves by the usual integral formula but they do not come from inner products in general. (Smoothly varying families of inner products in the tangent spaces of the base manifold constitute Riemannian metrics.)

The notion of generalized Berwald manifolds goes back to V. Wagner \cite{Wag1}. They are Finsler manifolds admitting linear connections such that the parallel transports preserve the Finslerian length of tangent vectors (compatibility condition). The first systematic investigation of generalized Berwald spaces is due to M. Hashiguchi \cite{H1} and \cite{H2}, M. Hashiguchi and Y. Ichijyo \cite{HY}, see also S. B\'{a}cs\'{o}, M. Hashiguchi and M. Matsumoto \cite{BHM}, S. Kikuchi \cite{Kikuchi} and M. Matsumoto \cite{M2}. They clarified that some special classes of generalized Berwald spaces can be given in terms of the conformal Finsler geometry (conformally Berwald and conformally flat Finsler spaces). The problem of the existence of the compatible linear connection and its intrinsic characterization persisted, however, in spite of the efforts being made by the classical approach. The new perspectives of the theory of generalized Berwald spaces appeared by introducing averaged Riemannian metrics and differential forms given by averaging. Using average quantities is a new and important trend in Finsler geometry \cite{V5}, \cite{V6}, \cite{Vin11} and \cite{V10}, R. G. Torrom\'{e} \cite{Torrome1}, T. Aikou \cite{Aikou1}, M. Crampin \cite{C2} and \cite{C1},  V. S. Matveev (et. al.: H-B. Rademacher, M. Troyanov) \cite{Mat2}, \cite{Mat1} and \cite{Mat3}. 
 
It is known \cite{V5} that any compatible linear connection must be metrical with respect to the averaged Riemannian metric given by integration of the Riemann--Finsler metric on the indicatrix hypersurfaces. Therefore any linear connection preserving the Finslerian length of tangent vectors is uniquely determined by its torsion. If the torsion is zero then we have a classical Berwald space and  the torsion-free compatible linear connection is the L\'{e}vi-Civita connection of the averaged Riemannian metric. It is a widely investigated area in Finsler geometry dominated by Z. I. Szab\'{o}'s famous structure theorem \cite{Szab1}. Otherwise, the torsion is a strange data we need to express in terms of intrinsic quantities. The problem of intrinsic characterization is also solved in the more general case of Finsler manifolds admitting semi-symmetric compatible linear connections \cite{V10}, see also \cite{V11}. The torsion tensor of a semi-symmetric compatible linear connection can be explicitly expressed in terms of the averaged Riemannian metric and differential forms given by averaging. Additionally, such a compatible linear connection is uniquely determined. 

The unicity of the compatible linear connection of a generalized Berwald manifold is false in general as the case of Randers metrics \cite{Vin1} and their generalizations show, see A. Tayebi and M. Barzegari \cite{Tayebi}. Difficulties coming from different possible solutions can be avoided by the idea of looking for the extremal solution in some sense: the extremal compatible linear connection of a generalized Berwald manifold keeps its torsion as close to zero as possible. (The bundle metric is induced by the averaged Riemannian metric on the base manifold.) Roughly speaking, we have to solve a hybrid conditional extremum problem at each point of the base manifold all of whose constraint equations involve functions defined on the indicatrix hypersurface (compatibility equations). The objective function to be minimized is a quadratic squared norm function defined on the finite dimensional fiber (vector space) of the torsion tensor bundle. Instead of the algorithmic solution \cite{V12} (based on the evaluation of the constraint equations at specially chosen elements of the tangent spaces point by point) the paper is devoted to the expression of the solution by using the method of Lagrange multipliers on function spaces point by point. We present a necessary and sufficient condition of the solvability and the solution is also given in terms of intrinsic quantities affecting the uniform size of the linear isometry groups of the indicatrices. If we have a generalized Berwald space, the pointwise solutions constitute a smooth section of the torsion tensor bundle. It is obviously a necessary and sufficient condition for a metric to be a generalized Berwald metric. We can weaken the  conditions by requiring the continuity of the section constituted by the pointwise solutions. A continuous section of the torsion tensor bundle implies continuous connection parameters on the base manifold. Solving the linear system of differential equations of the parallel vector fields, linear isometries can be constructed between the tangent spaces with respect to the Finsler metric (monochromatic metrics). According to  N. Bartelme{\ss}  and V. S. Matveev \cite{BM}, monochromatic metrics must be generalized Berwald metrics.  

Using the general method, we present explicit computations in case of spaces of dimension $2$ and $3$ including the special case of Randers metrics. They are motivated by the results of some previous investigations \cite{VO} and \cite{VO1}.

\section{Notations and terminology}

Let $M$ be a connected differentiable manifold with local coordinates $u^1, \ldots, u^n.$ The induced coordinate system of the tangent manifold $TM$ consists of the functions $x^1, \ldots, x^n$ and $y^1, \ldots, y^n$. For any $v\in T_pM$, $x^i(v):=u^i\circ \pi (v)=p^i$ and $y^i(v)=v(u^i)$, where $i=1, \ldots, n$ and $\pi \colon TM \to M$ is the canonical projection. 

A \emph{Finsler metric} is a non-negative continuous function $F\colon TM\to \mathbb{R}$ satisfying the following conditions:  $\displaystyle{F}$ is smooth on the complement of the zero section (regularity), $\displaystyle{F(tv)=tF(v)}$ for all $\displaystyle{t> 0}$ (positive homogeneity) and the Hessian 
$$g_{ij}=\frac{\partial^2 E}{\partial y^i \partial y^j}$$
of the energy function $E=F^2/2$ is positive definite at all nonzero elements $\displaystyle{v\in T_pM}$ (strong convexity), $p\in M$. The so-called \emph{Riemann--Finsler metric} $g$ is constituted by the components $g_{ij}$. It is defined on the complement of the zero section. The Riemann-Finsler metric makes the complement of the origin a Riemannian manifold in each tangent space. The canonical objects are the {\emph {volume form}} $\displaystyle{d\mu=\sqrt{\det g_{ij}}\ dy^1\wedge \ldots \wedge dy^n}$,
the \emph {Liouville vector field} $\displaystyle{C:=y^1\partial /\partial y^1 +\ldots +y^n\partial / \partial y^n}$ and the {\emph {induced volume form}}
$$\mu=\sqrt{\det g_{ij}}\ \sum_{i=1}^n (-1)^{i-1} \frac{y^i}{F} dy^1\wedge\ldots\wedge dy^{i-1}\wedge dy^{i+1}\ldots \wedge dy^n$$
 on the indicatrix hypersurface $\displaystyle{\partial K_p:=F^{-1}(1)\cap T_pM\ \  (p\in M)}$. The averaged Riemannian metric is defined by 
\begin{equation}
\label{averagemetric1}
\gamma_p (v,w):=\int_{\partial K_p} g(v, w)\, \mu=v^i w^j \int_{\partial K_p} g_{ij}\, \mu \ \ (v, w\in T_p M, p\in M).
\end{equation}

\begin{Def} A linear connection $\nabla$ on the base manifold $M$ is called \emph{compatible} with the Finsler metric if the parallel transports with respect to $\nabla$ preserve the Finslerian length of tangent vectors. Finsler manifolds admitting compatible linear connections are called generalized Berwald manifolds.
\end{Def}

Let $c\colon [0,1]\to M$ be a curve on the base manifold. Evaluating the Finsler metric along a parallel vector field $X$ with respect to $\nabla$ (a linear connection on the base manifold) we have that 
\begin{equation}
\label{eq:55}
(F \circ X)'={c^i}'\bigg(\frac{\partial F}{\partial x^i}-y^j {\Gamma}_{ij}^{k}\circ \pi \frac{\partial F}{\partial y^k}\bigg)\circ X.
\end{equation}
This means that the parallel transports with respect to $\nabla$ preserve the Finslerian length of tangent vectors (compatibility condition) if and only if
\begin{equation}
\label{cond1}
\frac{\partial F}{\partial x^i}-y^j {\Gamma}^k_{ij}\circ \pi \frac{\partial F}{\partial y^k}=0 \ \  (i=1, \ldots,n).
\end{equation}
The vector fields of type
\begin{equation}
\label{eq:6}
X_i^{h}:=\frac{\partial}{\partial x^i}-y^j {\Gamma}^k_{ij}\circ \pi \frac{\partial}{\partial y^k}
\end{equation}
span the horizontal distribution belonging to $\nabla$. We can introduce the volume forms $d\mu^*$, $\mu^*$ and the horizontal vector fields $X_i^{h^*}$ ($i=1, \ldots, n$) with respect to the L\'{e}vi-Civita connection of the averaged Riemannian metric $F^*(v):=\sqrt{\gamma_p(v,v)}$ in a similar way. 

\begin{Thm}
\label{heritage} \emph{\cite{V5}} If a linear connection on the base manifold is compatible with the Finsler metric function then it must be metrical with respect to the averaged Riemannian metric.
\end{Thm}

Using Theorem \ref{heritage} we can substitute the connection parameters with the components of the torsion tensor in the compatibility equations (\ref{cond1}). Following the standard Christoffel process 
\begin{equation}
\label{torsion}
\Gamma_{ij}^r=\Gamma_{ij}^{*r}-\frac{1}{2}\left(T^{l}_{jk}\gamma^{kr}\gamma_{il}+T^{l}_{ik}\gamma^{kr}\gamma_{jl}-T_{ij}^r\right)
\end{equation}
and the compatibility equations (\ref{cond1}) can be written into the form
\begin{equation}
\label{cond2}
X_i^{h^*}F+\frac{1}{2}y^j \left(T^{l}_{jk}\gamma^{kr}\gamma_{il}+T^{l}_{ik}\gamma^{kr}\gamma_{jl}-T_{ij}^r\right)\circ \pi \frac{\partial F}{\partial y^r}=0\ \  (i=1, \ldots,n).
\end{equation}
The correspondence $\nabla \rightleftharpoons T$ preserves the affine combinations of the linear connections, i.e. for any real number $\lambda\in \mathbb{R}$ we have 
$$\lambda \nabla_1+(1-\lambda)\nabla_2 \rightleftharpoons \lambda T_1+(1-\lambda)T_2.$$
Additionally, if $\nabla_1$ and $\nabla_2$ satisfy the compatibility equations (\ref{cond1}) then so does 
$$\lambda \nabla_1+(1-\lambda)\nabla_2.$$
This means that the torsion tensors of the compatible linear connections form an affine subspace $A_p$ in the linear space $\wedge^2 T_p^*M \otimes T_pM$ for any $p\in M$. As the point is varying we have an affine distribution of the torsion tensor bundle $\wedge^2 T^*M \otimes TM$. In terms of local coordinates, the bundle is spanned by
$$du^i \wedge du^j \otimes \frac{\partial}{\partial u^k} \quad (1\leq i < j\leq n, k=1, \ldots, n)$$
and its dimension is $\displaystyle{\binom{n}{2}n}$. 

Since the torsion tensor bundle can be equipped with a Riemannian metric in a natural way, we can measure the length of the torsion to formulate an extremum problem for the  compatible linear connection keeping its torsion as close to the origin as possible. 

\begin{Def} \label{bundle_metric} The products 
$$du^i \wedge du^j \otimes \frac{\partial}{\partial u^k} \quad (1\leq i< j \leq n, k=1, \ldots, n)$$
form an orthonormal basis at the point $p\in M$ if the coordinate vector fields $\partial/\partial u^1, \ldots, \partial/\partial u^n$ form an orthonormal basis with respect to the averaged Riemannian metric at the point $p\in M$. The norm of the torsion tensor is defined by
\begin{equation}
\label{torsionnorm}
\|T_p\|^2=\sum_{1\leq i < j \leq n} \sum_{k=1}^n {T_{ij}^k(p)}^2
\end{equation}
provided that 
$$T=\sum_{1\leq i < j \leq n} \sum_{k=1}^nT_{ij}^k du^i \wedge du^j \otimes \frac{\partial}{\partial u^k}$$ 
and the products form an orthonormal basis at the point $p\in M$. The corresponding inner product is
$$\langle T_p, S_p \rangle=\sum_{1\leq i < j \leq n} \sum_{k=1}^n T_{ij}^k(p)S_{ij}^k(p).$$
\end{Def}

In what follows we introduce the extremal compatible linear connection of a generalized Berwald manifold in terms of its torsion $T^0$ by taking the closest point $T^0_p \in A_p$ to the origin for any $p\in M$. 

\begin{Def} The extremal compatible linear connection of a generalized Berwald manifold is the uniquely determined compatible linear connection minimizing the norm of its torsion by taking the values of the pointwise minima. 
\end{Def}

\section{A conditional extremum problem for the extremal compatible linear connection} Let a point $p\in M$ be given. We are going to use the compatibility equations (\ref{cond2}) as constraint equations at the point $p\in M$. Therefore the conditional extremum problem
for the extremal compatible linear connection takes the form 
\begin{equation}
\label{condext} \min \frac{1}{2} \|T_p\|^2 \quad \textrm{subject to} \quad X_i^{h^*}F+\frac{1}{2}y^j \left(T^{l}_{jk}\gamma^{kr}\gamma_{il}+T^{l}_{ik}\gamma^{kr}\gamma_{jl}-T_{ij}^r\right)(p) \frac{\partial F}{\partial y^r}=0, 
\end{equation}
where $i=1, \ldots, n$. The coefficient of $T_{ab}^c$ ($1\leq a < b \leq n$, $c=1, \ldots, n$) is
\begin{equation}
\label{coeff}
\sigma_{c; i}^{ab}=\frac{1}{2}\left(\left(y^a \gamma^{br}-y^b\gamma^{ar}\right)\frac{\partial F}{\partial y^r}\gamma_{ic}+\left(\delta_i^a \gamma^{br}-\delta_i^b\gamma^{ar}\right)\frac{\partial F}{\partial y^r}y^j\gamma_{jc}-\left(\delta_i^a y^b-\delta_i^b y^a\right)\frac{\partial F}{\partial y^c}\right),
\end{equation}
where the index $i=1, \ldots, n$ refers to the corresponding compatibility equation. The symmetric differences in formula (\ref{coeff}) are due to $a< b$. If the coordinate vector fields form an orthonormal basis at $p\in M$ with respect to the averaged Riemannian metric $\gamma$, then 
\begin{equation}
\label{coeffort}
\sigma_{c; i}^{ab}=\frac{1}{2}\left(\delta_i^c\left(y^a \frac{\partial F}{\partial y^b}-y^b\frac{\partial F}{\partial y^a}\right)+\delta_i^a \left(y^c\frac{\partial F}{\partial y^b}-y^b\frac{\partial F}{\partial y^c}\right)-\delta_i^b \left(y^c\frac{\partial F}{\partial y^a}-y^a\frac{\partial F}{\partial y^c}\right)\right).
\end{equation}
The tensorial form of problem (\ref{condext}) is 
\begin{equation}
\label{condextLag}
\min \frac{1}{2} \|T_p\|^2 \quad \textrm{subject to} \quad X_i^{h^*}F+\langle T_p, \sigma_i \rangle=0 \quad (i=1, \ldots, n),
\end{equation}
where
$$\sigma_i=\sum_{1\leq a < b \leq n} \sum_{c=1}^n \sigma_{c; i}^{ab} dy^a\wedge dy^b \otimes \frac{\partial}{\partial y^c}$$
and the inner product on the vertical subspace over $\wedge^2 T^*_pM \otimes T_pM$ is induced by the Riemannian metric of the torsion tensor bundle:
$$\langle T_p, \sigma_i(v) \rangle :=\sum_{1\leq a < b \leq n} \sum_{c=1}^n T_{ab}^c(p) \sigma_{c; i}^{ab}(v);$$
note that the torsion tensor $T$ can be identified with its vertically lifted tensor in a canonical way: 
$$T=\sum_{1\leq a < b \leq n} \sum_{c=1}^n T_{ab}^c d u^a \wedge d u^b \otimes \frac{\partial}{\partial u^c} \rightleftharpoons \sum_{1\leq a < b \leq n} \sum_{c=1}^n T_{ab}^c\circ \pi  d y^a \wedge d y^b \otimes \frac{\partial}{\partial y^c}.$$

Using the homogeneity property of the Finsler metric, it is enough to consider the constraint equations over the Riemannian indicatrix hypersurface $S_p^*\subset T_pM$ with respect to the averaged Riemannian metric.

\subsection{The method of Lagrange multipliers} We introduce some concise notations to make the steps of the solution clear: 
\begin{itemize}
\item the fiber of the torsion tensor bundle at the point $p\in M$ is denoted by
$$V=\wedge^2 T^*_pM \otimes T_pM \sim \mathbb{R}^m.$$
It is a Euclidean vector space of dimension $\displaystyle{m=\binom{n}{2}n}$. The inner product is given by the averaged Riemannian metric in the sense of Definition \ref{bundle_metric}.
\item $S_p^* \subset T_pM \sim \mathbb{R}^n$ is the Riemannian unit sphere at the point $p\in M$ with respect to the averaged Riemannian metric. $\mathcal{C}^{\infty}(S_p^*)$ denotes the space of smooth functions defined on the Riemannian unit sphere. If necessary, we consider its completion as a metric space. The inner product is given by integration of the product of elements with respect to the induced volume form $\mu^*$ on the indicatrix hypersurface by the averaged Riemannian metric. 
\item The objective function is $f\colon V\to \mathbb{R}$, $f(T_p):=\|T_p\|^2/2$.
\item The constraints are  
$$g_i\colon V\to \mathcal{C}^{\infty}(S_p^*), \quad g_i(T_p):=X_i^{h^*}F+\langle T_p, \sigma_i \rangle \quad (i=1, \ldots, n)$$
and 
$$g\colon V\to W, \quad g(T_p):=\left(g_1(T_p), \ldots, g_n(T_p)\right),$$
where
$$W:=\mathcal{C}^{\infty}(S_p^*)\times \ldots \times \mathcal{C}^{\infty}(S_p^*) \quad (\textrm{$n$ times})$$
is a pre-Hilbert space equipped with the inner product
$$\langle \nu^1, \nu^2 \rangle=\sum_{i=1}^n \int_{S_p^*}\nu_i^1 \nu_i^2\, \mu^*.$$
If necessary, we consider its completion as a metric space.
\item $h^*:=(h_1^*, \ldots, h_n^*)=\left(X_1^{h^*}F, \ldots, X_n^{h^*}F\right)\in W$ (the affine term of the constraint).
\item $T_p^0$ is the solution of the conditional extremum problem
\begin{equation}
\label{condextLagnew}
\min f(T_p) \quad \textrm{subject to} \quad g(T_p)=0.
\end{equation}
\end{itemize}
Following the Lagrange principle we are looking for an element $\lambda$ of the dual space $W^*$ of $W$ such that
$$f'(T_p^0)=\lambda \circ g'(T_p^0).$$
Since
$$f'(T_p^0)(T_p)=\langle T_p, T_p^0\rangle$$
and
$$g'(T_p^0)(T_p)=\left(\langle T_p, \sigma_1 \rangle, \ldots, \langle T_p, \sigma_n \rangle \right)=\textrm{the linear part of $g$ at $T_p\in V$}$$
it follows, by the Riesz representation theorem, that
$$\langle T_p, T_p^0\rangle =\lambda\circ \left(\langle T_p, \sigma_1 \rangle, \ldots, \langle T_p, \sigma_n \rangle \right)=\int_{S_p^*}\lambda_1 \langle T_p, \sigma_1 \rangle\, \mu^*+ \ldots+\int_{S_p^*} \lambda_n \langle T_p, \sigma_n \rangle \, \mu^*=$$
$$\left \langle T_p, \int_{S_p^*}\lambda_1 \sigma_1 \, \mu^*\right \rangle + \ldots+ \left \langle T_p, \int_{S_p^*} \lambda_n\sigma_n \, \mu^* \right \rangle =\left \langle T_p, \int_{S_p^*}\lambda_1 \sigma_1\, \mu^* + \ldots+ \int_{S_p^*} \lambda_n\sigma_n \, \mu^*\right \rangle,$$
where the componentwise integration  
$$\int_{S_p^*}\lambda_i \sigma_i \, \mu^*= \sum_{1\leq a < b \leq n} \sum_{c=1}^n \int_{S_p^*}\lambda_i\sigma_{c; i}^{ab}\, \mu^* dy^a \wedge dy^b \otimes\frac{\partial}{\partial y^c}  \ \rightleftharpoons \ \sum_{1\leq a < b \leq n} \sum_{c=1}^n \int_{S_p^*}\lambda_i\sigma_{c; i}^{ab}\, \mu^* du^a \wedge du^b \otimes\frac{\partial}{\partial u^c}$$
allows us to consider the tensors on the base manifold by eliminating the formal dependence on the directional coordinates. Therefore 
\begin{equation}
\label{solutionform}T_p^0=T(\lambda):=\sum_{i=1}^n\int_{S_p^*}\lambda_i \sigma_i \mu^*.
\end{equation}

\subsection{Factorization and solvability} Using the constraint equations we have that for any $\mu=(\mu_1, \ldots, \mu_n)\in W$
$$\sum_{i=1}^n \mu_ig_i(T_p^0)=0 \quad \Leftrightarrow \quad \sum_{i=1}^n \mu_i h^*_i+\langle T_p^0, \sum_{i=1}^n \mu_i \sigma_i\rangle=0.$$
Integrating over $S_p^*$ it follows that
\begin{equation}\label{solvability01}
\langle \mu, h^*\rangle+\langle T_p^0, T(\mu)\rangle=0,
\end{equation}
where
\begin{equation}
\label{torsionmapping0} 
T\colon W\to V, \quad T(\mu)=\sum_{i=1}^n\int_{S_p^*}\mu_i \sigma_i \, \mu^*.
\end{equation}
Using (\ref{coeffort}), the components of the torsion tensor $T(\mu)$ are
\begin{equation}
\label{torsionmapping} T_{ab}^c(\mu)=\frac{1}{2}\int_{S_p^*} \mu_c\left(y^a \frac{\partial F}{\partial y^b}-y^b\frac{\partial F}{\partial y^a}\right)+\mu_a \left(y^c\frac{\partial F}{\partial y^b}-y^b\frac{\partial F}{\partial y^c}\right)-\mu_b \left(y^c\frac{\partial F}{\partial y^a}-y^a\frac{\partial F}{\partial y^c}\right)\, \mu^*
\end{equation}
$(1\leq a< b\leq n, c=1, \ldots, n)$. Factorizing $W$ by the kernel of $T$ we have an injective linear mapping
$$T\colon W /\textrm{\ Ker\ } T\to V, \quad T[\mu]=T(\mu),$$
where $V$ is a finite dimensional vector space. So is  $W /\textrm{\ Ker\ } T$. 

\begin{Thm}
\label{solvability03} The conditional extremum problem for the torsion of the extremal compatible linear connection is solvable at the point $p\in M$ if and only if the inclusion 
\begin{equation}
\label{solvability04}
\textrm{Ker\ } T\subset \ {h^*}^{\bot}
\end{equation}
is satisfied. In case of $h^*\neq 0$ and $\textrm{\ Ker\ } T\subset \ {h^*}^{\bot}$ we have that
$$\dim \textrm{Im\ }T =\dim W/\textrm{Ker\ } T= 1+\dim {h^*}^{\bot} /\textrm{\ Ker\ } T.$$
\end{Thm}

\begin{Pf} If the conditional extremum problem is solvable then, by substituting $\mu\in \textrm{Ker\ }T$ in (\ref{solvability01}), it follows that $\mu\in {h^*}^{\bot}$. If (\ref{solvability04}) holds, then we have a well-defined linear functional 
$$T(\mu) \mapsto - \left \langle \mu, h^* \right \rangle $$
on the image $\textrm{Im\ }T$. Using the Riesz representation theorem, there is a uniquely determined element $T_p^0=T(\lambda)\in \textrm{Im\ }T$ such that
$$-\langle \mu, h^*\rangle=\langle T_p^0, T(\mu)\rangle \ \Leftrightarrow \ \langle \mu, h^*\rangle+\langle T_p^0, T(\mu)\rangle=0 \ \Leftrightarrow \ \sum_{i=1}^n \int_{S_p^*}\mu_ig_i(T_p^0)\, \mu^*=0$$
for any $\mu\in W$. Therefore $g(T_p^0)=0$ and we have a non-empty solution space. In case of $h^*\neq 0$ and $\textrm{Ker\ } T\subset \ {h^*}^{\bot}$, the factor space  ${h^*}^{\bot} /\textrm{\ Ker\ } T$ is a hyperspace in $W /\textrm{\ Ker\ } T\cong \textrm{\ Im\ }T$. Indeed, (\ref{solvability01}) implies that the image of ${h^*}^{\bot}$ is the orthogonal complement of $T_p^0=T(\lambda)\neq 0$ in $\textrm{Im\ } T$ and the factorization provides the one-to-one correspondence. 
\end{Pf}

\begin{Rem}\label{riemann}{\emph{The torsion mapping \eqref{torsionmapping0} is identically zero at the point $p\in M$ if and only if the Finsler  indicatrix reduces to a quadratic hypersurface. Indeed, by formula \eqref{torsionmapping}, the vanishing of the torsion mapping at the point $p\in M$ implies that  
\begin{equation}
\label{contracting}
y^a \frac{\partial F}{\partial y^b}-y^b\frac{\partial F}{\partial y^a}=0 
\end{equation}
for any $a$, $b=1, \ldots, n$. Contracting equation (\ref{contracting}) 
$$y^a =\sum_{b=1}^n \left(y^b\right)^2\frac{\partial \log F}{\partial y^a}\ \ \Leftrightarrow \ \ \frac{y^a}{\sum_{b=1}^n \left(y^b\right)^2}=\frac{\partial \log F}{\partial y^a}\ \ \Leftrightarrow \ \ \frac{\partial \log F^*}{\partial y^a}=\frac{\partial \log F}{\partial y^a}\ \ \Leftrightarrow \ \ F=e^{\alpha(p)}F^*,$$
i.e. $T_pM$ is a vertically contact tangent space \cite{V12}: the Finsler indicatrix is homothetic to a Riemannian one at the point $p\in M$. Having parallel transports with respect to a compatible linear connection (generalized Berwald metrics), such a geometric property can be extended to the entire (connected) base manifold. Therefore the generalized Berwald Finsler metric reduces to a Riemannian metric.}}
\end{Rem}

\begin{Rem}{\emph{Assuming the solvability condition (\ref{solvability04}) and $T_p^0=T(\lambda)$, the Cauchy--Bunyakovsky--Schwarz inequality implies  that 
$$\langle \mu, h^* \rangle^2\stackrel{(\ref{solvability01})}{=}\langle T(\lambda), T(\mu)\rangle^2\leq \|T(\lambda)\|^2 \ \|T(\mu)\|^2.$$
In what follows we suppose that the torsion mapping \eqref{torsionmapping0} is not identically zero to avoid the trivial case of the   quadratic indicatrices (see Remark \ref{riemann}). Since we have equality under the choice of $\lambda=\mu$, $\lambda=(\lambda_1, \ldots, \lambda_n)$ is the position where the supremum
\begin{equation}
\label{support01}
\|T(\lambda)\|^2=\sup_{\|T(\mu)\|^2\neq 0}\langle \mu, h^*\rangle^2/\|T(\mu)\|^2
\end{equation} 
is attained at. We can also prove that $\displaystyle{D:=\{\mu \in W\ | \ \|T(\mu)\|^2\leq 1\}}$ is a convex subset in $W$ containing the origin in its interior because of
$$\|T(\mu)\|^2= \| \sum_{i=1}^n\int_{S_p^*}\mu_i \sigma_i \, \mu^* \|^2\leq \|\mu\|^2 \sum_{i=1}^n \int_{S_p^*}\|\sigma_i\|^2\, \mu^*,$$
i.e. $\displaystyle{\|T(\mu)\|^2<1}$ provided\footnote{Note that 
$$\sum_{i=1}^n \int_{S_p^*}\|\sigma_i\|^2\, \mu^*=0$$
is equivalent to the trivial case of the quadratic indicatrices (see Remark \ref{riemann}).} that
$$\sum_{i=1}^n \int_{S_p^*}\|\sigma_i\|^2\, \mu^*\neq 0 \quad \textrm{and} \quad \|\mu\|^2 < \frac{1}{\sum_{i=1}^n \int_{S_p^*}\|\sigma_i\|^2\, \mu^*}.$$
Therefore 
\begin{equation}
\label{support02}
\|T(\lambda)\|^2=\sup_{\|T(\mu)\|^2\neq 0}\langle \mu, h^*\rangle^2/\|T(\mu)\|^2 = \sup_{\|T(\mu)\|^2=1}\langle \mu, h^*\rangle^2=L_{D^*}^2(h^*),
\end{equation}
where $L_{D^*}$ is the Minkowski functional of the polar set $D^*$.}}
\end{Rem}

\subsection{Formal expressions of the solution} Assuming the solvability condition (\ref{solvability04}) we can formulate a finite dimensional problem by taking
$$T_p^0\stackrel{(\ref{solutionform})}{=}T(\lambda)=\sum_{j=1}^k r_jT(\mu^j), \quad \textrm{where}\quad \dim W /\textrm{\ Ker\ } T=k,$$
$[\mu^1], \ldots, [\mu^k]$ is a basis of $W /\textrm{\ Ker\ } T$ and $r_1, \ldots, r_k$ are real numbers. Substituting  $T(\lambda)$ into the constraint equations of the form (\ref{solvability01}) we have that
$$\langle \mu^1, h^*\rangle+\sum_{j=1}^k r_j\langle T(\mu^j), T(\mu^1)\rangle=0, \ \ldots, \ \langle \mu^k, h^*\rangle +\sum_{j=1}^k r_j\langle T(\mu^j), T(\mu^k)\rangle=0$$
and the uniquely determined solutions $r_1, \ldots, r_k$ can be expressed in terms of the inverse Gram-matrix of $T(\mu^1), \ldots, T(\mu^k)$. The "orthogonalization process"
\begin{equation}
\label{ortproc}
\nu^1:=\mu^1, \ \nu^2=\mu^2-\frac{\langle T(\mu_2), T(\nu_1)\rangle}{\|T(\nu_1)\|^2}\nu_1, \ \ldots, \ \nu^k=\mu^k-\sum_{j=1}^{k-1}\frac{\langle T(\mu_k), T(\nu_j)\rangle}{\|T(\nu_j)\|^2}\nu_j
\end{equation}
gives an orthogonal basis $T(\nu^1), \ldots, T(\nu^k)$ in $\textrm{\ Im}\ T$. Taking 
$$T_p^0\stackrel{(\ref{solutionform})}{=}T(\lambda)=\sum_{j=1}^k s_jT(\nu^j)  \quad \textrm{where}\quad \dim W /\textrm{\ Ker\ } T=k $$
and $s_1, \dots, s_k$ are real numbers. The constraint equations of the form (\ref{solvability01}) imply that
\begin{equation}
\label{solution}
T(\lambda)=-\sum_{i=1}^k \frac{\langle \nu^i, h^* \rangle}{\|T(\nu^i)\|^2}T(\nu^i).
\end{equation} 

\subsection{Computations in 2D} In case of two-dimensional spaces
$$T(\mu)=\sum_{i=1}^2\int_{S_p^*}\mu_i \sigma_i \, \mu^*=\sum_{1\leq a < b \leq 2} \sum_{c=1}^2 \sum_{i=1}^2\int_{S_p^*}\mu_i\sigma_{c; i}^{ab}\, \mu^* du^a \wedge du^b \otimes\frac{\partial}{\partial u^c}=$$
$$\sum_{c=1}^2 \sum_{i=1}^2\int_{S_p^*}\mu_i \sigma_{c; i}^{12}\, \mu^*du^1\wedge du^2\otimes \frac{\partial}{\partial u^c}=$$
$$\sum_{i=1}^2\int_{S_p^*}\mu_i \sigma_{1; i}^{12}\, \mu^*du^1\wedge du^2\otimes \frac{\partial}{\partial u^1}+\sum_{i=1}^2\int_{S_p^*}\mu_i \sigma_{2; i}^{12}\, \mu^*du^1\wedge du^2\otimes \frac{\partial}{\partial u^2}=$$
$$\int_{S_p^*}\mu_1 \sigma_{1; 1}^{12}+\mu_2 \sigma_{1; 2}^{12} \, \mu^*du^1\wedge du^2\otimes \frac{\partial}{\partial u^1}+\int_{S_p^*}\mu_1 \sigma_{2; 1}^{12}+\mu_2 \sigma_{2; 2}^{12}\, \mu^*du^1\wedge du^2\otimes \frac{\partial}{\partial u^2}\stackrel{(\ref{coeffort})}{=}$$
$$\int_{S_p^*}\mu_1 f_{12}\, \mu^*du^1\wedge du^2\otimes \frac{\partial}{\partial u^1}+\int_{S_p^*}\mu_2 f_{12}\, \mu^*du^1\wedge du^2\otimes \frac{\partial}{\partial u^2},$$
where
$$f_{12}=y^1 \frac{\partial F}{\partial y^2}-y^2\frac{\partial F}{\partial y^1}.$$
According to the expression of $T(\mu)$ let us choose the elements $\nu^1=(f_{12}, 0)$ and $\nu^2=(0, f_{12})$ to provide the orthogonal basis
$$T(\nu^1)=\int_{S_p^*}\left(f_{12}\right)^2 \, \mu^* du^1\wedge du^2\otimes \frac{\partial}{\partial u^1},\ T(\nu^2)=\int_{S_p^*}\left(f_{12}\right)^2 \, \mu^* du^1\wedge du^2\otimes \frac{\partial}{\partial u^2}$$
in the fiber of the torsion tensor bundle at the point $p\in M$. Using formula (\ref{solution}) we have
\begin{equation}
\label{2dsolution}
T(\lambda)=-\frac{1}{\int_{S_p^*}\left(f_{12}\right)^2\, \mu^*}\left(\int_{S_p^*} h_1^*f_{12}\, \mu^*du^1\wedge du^2\otimes \frac{\partial}{\partial u^1}+ \int_{S_p^*} h_2^*f_{12}\, \mu^*du^1\wedge du^2\otimes \frac{\partial}{\partial u^2}\right).
\end{equation}
The extremal compatible linear connection is the only possible solution of the compatibility equations because
$$\dim W/ \textrm{\ Ker} \ T =\dim \textrm{Im\ }T=\dim V=2$$
provided that the denominator $\int_{S_p^*}\left(f_{12}\right)^2\, \mu^*\neq 0$. If the denominator is zero then we have the trivial case of the quadratic indicatrices (see Remark \ref{riemann}). 

\begin{Thm}
\label{2d} The torsion of the extremal compatible linear connection of a connected non-Riemannian generalized Berwald surface is  
$$-\frac{1}{\int_{S_p^*}\left(f_{12}\right)^2\, \mu^*}\left(\int_{S_p^*} h_1^* f_{12}\, \mu^* du^1\wedge du^2\otimes \frac{\partial}{\partial u^1}+ \int_{S_p^*} h_2^* f_{12}\, \mu^* du^1\wedge du^2\otimes \frac{\partial}{\partial u^2}\right) \quad (p\in M)$$
where
$$h_1^*=X_1^*F ,\ h_2^*=X_2^*F, \ f_{12}=y^1 \frac{\partial F}{\partial y^2}-y^2\frac{\partial F}{\partial y^1}$$
and the coordinate vector fields form an orthonormal system with respect to the averaged Riemannian metric at the point $p\in M$. 
\end{Thm}            

Under the condition
\begin{equation}
\label{allornot2d}
\int_{S_p^*}\left(f_{12}\right)^2\, \mu^*\neq 0 \quad (p\in M),
\end{equation}
the (pointwise) extremal solutions form a smooth section of the torsion tensor bundle. If the corresponding linear connection is compatible with the Finsler metric then we have a non-Riemannian generalized Berwald surface: recall that condition
$$ \int_{S_p^*}\left(f_{12}\right)^2\, \mu^* = 0$$
holds either everywhere (Riemannian surfaces) or nowhere (see Remark \ref{riemann}) for a generalized Berwald metric.

\subsection{The main results: explicit expressions of the solution} \label{expsol} We are going to investigate the general case motivated by the computations in case of dimension two. Let us introduce the quantities 
$$f_{ab}=y^a \frac{\partial F}{\partial y^b}-y^b\frac{\partial F}{\partial y^a}\in \mathcal{C}^{\infty} (S_p^*) \quad (1\leq a < b \leq n)$$
at the point $p\in M$ with respect to the lexicographic ordering
\begin{equation}
\label{system01}
f_{12}, \ldots, f_{1n}, f_{23}, \ldots, f_{2n}, \ldots, f_{n-1 n}.
\end{equation}
Its Gram matrix is
$$\mathcal{G}(f_{ab})=\left(\int_{S_p^*} f_{ab}f_{kl} \, \mu^*\right)_{i=(a-1)n-\frac{(a-1)a}{2}+b-a \  j=(k-1)n-\frac{(k-1)k}{2}+l-k}\ \in \textrm{\ Mat}_{\binom{n}{2}\times \binom{n}{2}}(\mathbb{R}),$$
where the correcting terms $\displaystyle{\frac{(a-1)a}{2}}$ and $\displaystyle{\frac{(k-1)k}{2}}$ are due to $a<b$ and $k<l$, respectively:
$$(a-1)n-\frac{(a-1)a}{2}=(a-1)n-\left(1+\ldots+(a-1)\right),$$
$$(k-1)n-\frac{(k-1)k}{2}=(k-1)n-\left(1+\ldots+(k-1)\right).$$
The corresponding system of elements in $W$ is
\begin{equation}\label{system02}
\begin{aligned}
\mu_{12}^1=(f_{12}, 0, \ldots, 0), \ \mu_{13}^1=(f_{13}, 0, \ldots, 0), \ldots,\ \mu_{1n}^1=(f_{1n}, 0, \ldots, 0),\\
\mu_{23}^1=(f_{23}, 0, \ldots, 0), \ \mu_{24}^1=(f_{24}, 0, \ldots, 0), \ldots,\ \mu_{n-1 n}^1=(f_{n-1 n}, 0, \ldots, 0),\\
\mu_{12}^2=(0, f_{12}, 0, \ldots, 0), \ \mu_{13}^2=(0, f_{13}, 0, \ldots, 0), \ldots,\ \mu_{1n}^2=(0, f_{1n}, 0, \ldots, 0), \ldots,
\end{aligned}
\end{equation}
i.e. $\mu_{kl}^j$ is an element in $W$ with zero coordinates except $f_{kl}$ in the $j$-th position. The sequence number of the element according to the lexicographic ordering is 
$$(j-1) \binom{n}{2}+(k-1)n-\frac{(k-1)k}{2}+l-k.$$
The Gram matrix of the system is a block-diagonal matrix in $\textrm{\ Mat}_{n\binom{n}{2}\times n\binom{n}{2}}(\mathbb{R})$ with the Gram matrix $\mathcal{G}(f_{ab})$ along the diagonal ($n$ times). Therefore  
\begin{equation}
\label{Gramsrel}
\textrm{\ \emph{rank}}\ \mathcal{G}(\mu_{kl}^j)=n \cdot \textrm{\emph{rank}} \ \mathcal{G}(f_{ab}).
\end{equation}
Especially,
$$\left |\mathcal{G}(\mu_{kl}^j)\right|=\left|\mathcal{G}(f_{ab})\right|^n.$$

The next result clarifies how the rank of the Gram matrix $\mathcal{G}(f_{ab})$ can decrease.

\begin{Thm}
\label{linisom}
The rank of the Gram matrix $\mathcal{G}(f_{ab})$ is
$$\textrm{\ \emph{rank}}\ \mathcal{G}(f_{ab})=\binom{n}{2}-d,$$
where $d$ is the dimension of the linear isometry group of the Finslerian indicatrix at the point $p\in M$.
\end{Thm}

\begin{Pf}
The linear isometry group of the Finslerian indicatrix at the point $p\in M$ consists of linear transformations $\varphi\colon T_pM\to T_pM$ leaving the Finslerian indicatrix invariant. In other words
$$\varphi^* F:=F\circ \varphi=F \quad \Leftrightarrow \quad \varphi^* g=g,$$
because $\varphi$ is linear and the Riemann--Finsler metric is constituted by the second order directional derivatives of the energy function $E:=F^2/2$ satisfying $\varphi^* E:=E\circ \varphi=E.$ Therefore $\varphi$ is an orthogonal transformation with respect to the averaged Riemannian metric at the point $p\in M$ and the linear isometry group of the Finslerian indicatrix is a compact subgroup of the orthogonal group. Using the closed subgroup theorem, it is a compact Lie group of dimension $d$. To finish the proof we are going to give the geometric interpretation of equation
$$\sum_{1\leq k< l\leq n} r_{kl}f_{kl}=0,$$
where the left hand side is a non-trivial linear combination of the system (\ref{system01}) with real scalars. Since
$$\sum_{1\leq k< l\leq n} r_{kl}f_{kl}=\sum_{j=1}^n \left(\sum_{1\leq k< j}r_{kj}y^k-\sum_{j < l\leq n}r_{jl}y^l\right)\frac{\partial F}{\partial y^j}=
\left(0-\sum_{1 < l\leq n} r_{1l}y^l\right)\frac{\partial F}{\partial y^1}+$$
$$\left(\sum_{1\leq k< 2}r_{k2}y^k-\sum_{2 < l\leq n}r_{2l}y^l\right)\frac{\partial F}{\partial y^2}+\ldots +\left(\sum_{1\leq k< n} r_{kn}y^k-0\right)\frac{\partial F}{\partial y^n},$$
the vanishing of the linear combination can be written as
\begin{equation}
\label{isom01}
\gamma_p\left(A\circ y, \textrm{\ grad} \ F\right)=0,
\end{equation}
where $A$ is a skew-symmetric matrix,
$$y=(y^1, \ldots, y^n) \quad \textrm{and} \quad \textrm{\ grad} \ F=\left(\frac{\partial F}{\partial y^1}, \ \ldots, \ \frac{\partial F}{\partial y^n}\right)$$
is the Euclidean gradient of the Finsler metric $F$ with respect to the averaged Riemannian metric at the point $p\in M$. Equation \eqref{isom01} means that $F$ is constant along the integral curve defined by $A \circ c(t)=c'(t)$, i.e. the one-parameter subgroup $\varphi_t:=\exp (tA)$ leaves the Finslerian indicatrix $S_p=F^{-1}(1)\cap T_pM$ invariant. Therefore vanishing linear combinations correspond to the elements of the Lie algebra of the linear isometry group of the Finslerian indicatrix and vice versa. 
\end{Pf}

\begin{Cor}
\label{cor:01}
For a non-Riemannian generalized Berwald manifold the Gram matrix $\mathcal{G}(f_{ab})$ is of constant rank between
$$n-1\leq \textrm{\ \emph{rank}}\ \mathcal{G}(f_{ab})\leq\binom{n}{2}\ \ \Leftrightarrow \ \ n(n-1)\leq \textrm{\ \emph{rank}}\ \mathcal{G}(\mu_{kl}^j)\leq n\binom{n}{2}.$$
If the rank is maximal, then any compatible linear connection is of zero curvature.
\end{Cor}

\begin{Pf}
The rank is constant because the linear isometry groups of the Finslerian indicatrices at different points of a connected base manifold  are conjugated by the parallel transport induced by a compatible linear connection. Using Wang's theorem 
\cite{V11a} (see Remark 7) we have that the linear isometry group of the Finslerian indicatrix is of dimension at most $(n-1)(n-2)/2$ unless it is a quadratic hypersurface (Riemannian case). Finally, the holonomy group of any compatible linear connection is a subgroup of the linear isometry group of the Finslerian indicatrix at any point of the base manifold. If we have a trivial linear isometry group then the holonomy group is also trivial and the corresponding compatible linear connection is of zero curvature because of the holonomy theorem due to W. Ambrose and I. M. Singer.
\end{Pf}

\begin{Thm}
\label{nd} The torsion of the extremal compatible linear connection of a connected non-Riemannian generalized Berwald space is of the form
$$T(\lambda)=\sum_{j}\sum_{k < l} r_{kl}^jT(\mu_{kl}^j),$$
where the summation is taken with respect to the indices of a maximal linearly independent system of elements $\mu_{kl}^j$ $(1\leq k < l \leq n, \ j=1,\ 2, \ \ldots, \ n)$ and the coordinate vector fields form an orthonormal system with respect to the averaged Riemannian metric at the point $p\in M$. Especially, the coefficients $r_{kl}^j$ can be expressed in terms of the Gram matrix of the maximal linearly independent system of elements $\mu_{kl}^j$ $(1\leq k < l \leq n, \ j=1,\ \ldots, \ n)$ and the inner product of the elements $\mu_{kl}^j$ with $h^*$:
$$\begin{bmatrix}\vdots\\[2pt]r_{kl}^j\\[2pt] \vdots\end{bmatrix}=-\mathcal{G}^{-1}\left(T(\mu_{kl}^j)\right)\begin{bmatrix}\vdots\\[2pt] \langle \mu_{kl}^j, h^*\rangle \\[2pt] \vdots\end{bmatrix}.$$
\end{Thm}

For the proof see the next section. Under the condition
\begin{equation}
\label{allornotnd}
\textrm{rank\ } \mathcal{G}(f_{ab})=\textrm{const.\ }\geq n-1 \ (p\in M) \ \Leftrightarrow\ \ 
\textrm{rank\ } \mathcal{G}(\mu_{kl}^j)\stackrel{\eqref{Gramsrel}}{=}\textrm{const.\ }\geq n(n-1) \ (p\in M),
\end{equation}
the (pointwise) extremal solutions form a smooth section of the torsion tensor bundle. If the corresponding linear connection is compatible with the Finsler metric then we have a non-Riemannian generalized Berwald space: recall that condition 
$$\textrm{rank\ } \mathcal{G}(f_{ab})=\textrm{const.\ } \ (p\in M) \ \Leftrightarrow\ \ \textrm{rank\ } \mathcal{G}(\mu_{kl}^j)\stackrel{\eqref{Gramsrel}}{=}\textrm{const.\ } \ (p\in M)$$
obviously holds for a generalized Berwald metric because the parallel transports with respect to a compatible linear connection provide linear isometries between the tangent spaces. If the lower bound in \eqref{allornotnd} does not work then we have a Riemannian space because of Wang's theorem for the maximal dimension of the linear isometry group of a non-Riemannian indicatrix hypersurface (see Corollary \ref{cor:01}):
$$d\leq \frac{(n-1)(n-2)}{2}\ \ \Rightarrow\ \ \textrm{rank\ } \mathcal{G}(f_{ab})=\textrm{const.\ }=\binom{n}{2}-d\geq n-1.$$

Using the orthogonalization process (\ref{ortproc}) we can write that 
$$T(\lambda)=-\sum_{j}\sum_{k < l} \frac{\langle \nu_{kl}^j, h^*\rangle}{\|T(\nu_{kl}^j)\|^2}T(\nu_{kl}^j).$$
Especially, if the rank of the system (\ref{system02}) is maximal then the extremal compatible linear connection is the only compatible linear connection because of
$$\dim W/\textrm{\ Ker}\ T=\dim \textrm{\ Im\ } T=n\binom{n}{2}=\dim V.$$
According to Corollary \ref{cor:01}, it is of zero curvature and its torsion is of the form 
$$T(\lambda)=-\sum_{j=1}^n\sum_{1\leq k < l \leq n} \frac{\langle \nu_{kl}^j, h^*\rangle}{\|T(\nu_{kl}^j)\|^2}T(\nu_{kl}^j).$$ 

\section{The proof of Theorem \ref{nd}}

In what follows we are going to investigate the rank of the system $T(\mu_{kl}^j)$ to prove that it is closely related to the rank of the Gram matrix $\mathcal{G}(f_{ab})$. In case of maximal rank the extremal compatible linear connection is the only possible solution of the compatibility equations. Additional solutions are possible if the dimension of the space $W/ \textrm{\ Ker}\ T$ decreases. It corresponds to the decreasing rank of $\mathcal{G}(f_{ab})$ due to the appearance of linear isometries of the Finslerian indicatrix at the point $p\in M$. 

\begin{Thm} \label{dimdir}
The direction space given by the constraint equations is of dimension at least $nd$, where $d$ is the dimension of the linear isometry group of the Finslerian indicatrix at the point $p\in M$.
\end{Thm}

\begin{Pf} Using \eqref{condext} the direction space is given by
$$\frac{1}{2}y^j \left(T^{l}_{jk}\gamma^{kr}\gamma_{il}+T^{l}_{ik}\gamma^{kr}\gamma_{jl}-T_{ij}^r\right)(p) \frac{\partial F}{\partial y^r}=0 \quad (i=1, \ldots, n).$$
In terms of an orthonormal system of the coordinate vector fields at $p\in M$ we have that
\begin{equation}
\label{dimdir01}
\frac{1}{2} y^j \left(T^{i}_{jr}+T^{l}_{ir}\delta_{jl}-T_{ij}^r\right)(p) \frac{\partial F}{\partial y^r}=0 \quad (i=1, \ldots, n).
\end{equation}
Comparing equations \eqref{isom01} and \eqref{dimdir01}, equality 
$$\frac{1}{2}\left(T^{i}_{jr}+T^{j}_{ir}-T_{ij}^r\right)(p)=A_{jr}$$
provides a correspondence between the elements of the direction space given by the constraint equations and the elements of the Lie algebra of the linear isometry group. Since both $i=1, \ldots, n$ and $A$ can be independently chosen, it follows that the dimension of the direction space is at least $nd$.
\end{Pf}

\begin{Cor} 
\label{cor:02} The solution of the conditional extremum problem is in $\textrm{\ Im\ } T$ of dimension at most
$$n\left(\binom{n}{2}-d\right).$$
\end{Cor}

\begin{Pf}
The constraint equations of the form (\ref{solvability01}) show that the direction space of the torsion tensors of the compatible linear connections at $p\in M$ is defined by
$$\langle T, T(\mu)\rangle=0 \ \ (\mu\in W)$$
and, consequently, it is the orthogonal complement of $\textrm{\ Im\ } T$. 
\end{Pf}

Using elements of the form $T(\mu_{kl}^j)$ for the explicit expression of the solution $T_p^0=T(\lambda)$ is motivated by \eqref{Gramsrel}, Theorem \ref{linisom} and Corollary \ref{cor:02}.

\subsection{Computations in 3D} By the help of formula \eqref{torsionmapping}, direct computations give the following table (integration is taken on the Riemannian indicatrix with respect to the induced volume form at the point $p\in M$). The rank of the matrix table is equal to the rank of the Gram matrix of system (\ref{system02}) because they can be  transformed into each other by performing rank preserving operations on the rows (addition, subtraction and exchange of rows): adding/subtracting row $7$ to row $3$/from row $5$, adding the (new) fifth and third rows with alternating signs to two times the seventh row, exchanging to put the new rows $3$, $5$ and $7$ in the corresponding positions of the  Gram matrix.

\begin{center}
{\tiny \addtolength{\tabcolsep}{-2pt}{
\begin{tabular}{|c|c|c|c||c|c|c||c|c|c|}
\hline
&&&&&&&&&\\                                                                
$T_{ab}^c$&$T(\mu_{12}^1)$&$T(\mu_{13}^1)$&$T(\mu_{23}^1)$&$T(\mu_{12}^2)$&$T(\mu_{13}^2)$&$T(\mu_{23}^2)$&$T(\mu_{12}^3)$&$T(\mu_{13}^3)$&$T(\mu_{23}^3)$\\
&&&&&&&&&\\
\hline
&&&&&&&&&\\
$T_{12}^1$ & $\int(f_{12})^2 $&$\int f_{12}f_{13}$&$\int f_{12}f_{23}$&0&0&0&0&0& 0\\
&&&&&&&&&\\
\hline 
&&&&&&&&&\\
$T_{13}^1$ &$\int f_{12}f_{13}$ &$\int(f_{13})^2$&$\int f_{13}f_{23}$&0&0&0&0&0&0\\
&&&&&&&&&\\
\hline
&&&&&&&&&\\
$T_{23}^1$ &$\int f_{12}f_{23}/2$ &$\int f_{13}f_{23}/2$&$\int (f_{23})^2/2$&$\int f_{12}f_{13}/2$&$\int (f_{13})^2/2$&$\int f_{13}f_{23}/2$&$-\int(f_{12})^2/2$&$-\int f_{12}f_{13}/2$&$-\int f_{12}f_{23}/2$\\
&&&&&&&&&\\
\hline
&&&&&&&&&\\
$T_{12}^2$ &0 &0&0&$\int(f_{12})^2 $&$\int f_{12}f_{13}$&$\int f_{12}f_{23}$&0&0&0\\
&&&&&&&&&\\
\hline
&&&&&&&&&\\
$T_{13}^2$ &$\int f_{12}f_{23}/2$ &$\int f_{13}f_{23}/2$&$\int (f_{23})^2/2$&$\int f_{12}f_{13}/2$&$\int (f_{13})^2/2$&$\int f_{13}f_{23}/2$&$\int(f_{12})^2/2$&$\int f_{12}f_{13}/2$&$\int f_{12}f_{23}/2$\\
&&&&&&&&&\\
\hline
&&&&&&&&&\\
$T_{23}^2$ &0 &0&0&$\int f_{12}f_{23}$ &$\int f_{13}f_{23}$&$\int(f_{23})^2$&0&0&0\\
&&&&&&&&&\\
\hline
&&&&&&&&&\\
$T_{12}^3$ &$-\int f_{12}f_{23}/2$ &$-\int f_{13}f_{23}/2$&$-\int (f_{23})^2/2$&$\int f_{12}f_{13}/2$&$\int (f_{13})^2/2$&$\int f_{13}f_{23}/2$&$\int(f_{12})^2/2$&$\int f_{12}f_{13}/2$&$\int f_{12}f_{23}/2$\\
&&&&&&&&&\\
\hline
&&&&&&&&&\\
$T_{13}^3$ &0 &0&0&0&0&0&$\int f_{12}f_{13}$ &$\int(f_{13})^2$&$\int f_{13}f_{23}$\\
&&&&&&&&&\\
\hline
&&&&&&&&&\\
$T_{23}^3$ &0 &0&0&0&0&0&$\int f_{12}f_{23}$ &$\int f_{13}f_{23}$&$\int(f_{23})^2$\\
&&&&&&&&&\\
\hline
\end{tabular}}}
\end{center}

Another way to prove rank equivalence is to show correspondence between the vanishing non-trivial linear combinations of the rows. The rows of the matrix table are constituted by the vectors 
$$
a_1=\left(\int(f_{12})^2, \int f_{12}f_{13}, \int f_{12}f_{23}\right),\ a_2=\left(\int f_{12}f_{13}, \int(f_{13})^2, \int f_{13}f_{23}\right),$$
$$
a_3=\left(\int f_{12}f_{23}, \int f_{13}f_{23}, \int(f_{23})^2\right).
$$
Therefore a vanishing non-trivial linear combination of the rows of the matrix table can be written as 
$$\lambda_1 a_1+\lambda_2 a_2+\frac{\lambda_3+\lambda_5-\lambda_7}{2}a_3=(0,0,0), \ \lambda_4 a_1+\frac{\lambda_3+\lambda_5+\lambda_7}{2}a_2+\lambda_6 a_3=(0,0,0),$$
$$\frac{-\lambda_3+\lambda_5+\lambda_7}{2}a_1+\lambda_8 a_2+\lambda_9 a_3=(0,0,0)
$$
corresponding to a non-trivial vanishing linear combination of the rows of the Gram matrix of system \eqref{system02} (it is a block-diagonal matrix with the Gram matrix $\mathcal{G}(f_{ab})$ along the diagonal $3$ times) whose coefficients are
$$\lambda_1, \ \lambda_2, \ \lambda_3'=\frac{\lambda_3+\lambda_5-\lambda_7}{2},\ \lambda_4, \ \lambda_5'=\frac{\lambda_3+\lambda_5+\lambda_7}{2}, \ \lambda_6,\ \lambda_7'=\frac{-\lambda_3+\lambda_5+\lambda_7}{2}, \ \lambda_8, \ \lambda_9.$$

\begin{Thm}
\label{3d} The torsion of the extremal compatible linear connection of a connected non-Riemannian generalized Berwald space of dimension three  is of the form
$$T(\lambda)=\sum_{j}\sum_{k < l} r_{kl}^jT(\mu_{kl}^j),$$
where the summation is taken with respect to the indices of a maximal linearly independent system of elements $\mu_{kl}^j$ $(1\leq k < l \leq 3, \ j=1,\ 2, \ 3)$ and the coordinate vector fields form an orthonormal system with respect to the averaged Riemannian metric at the point $p\in M$. Especially, the coefficients $r_{kl}^j$ can be expressed in terms of the Gram matrix of the maximal linearly independent system of elements $\mu_{kl}^j$ $(1\leq k < l \leq 3, \ j=1,\ 2, \ 3)$ and the inner product of the elements $\mu_{kl}^j$ with $h^*$:
$$\begin{bmatrix}\vdots \\[2pt] r_{kl}^j \\[2pt] \vdots \end{bmatrix}=-\mathcal{G}^{-1}\left(T(\mu_{kl}^j)\right)\begin{bmatrix}\vdots \\[2pt] \langle \mu_{kl}^j, h^*\rangle \\[2pt] \vdots \end{bmatrix}.$$
\end{Thm}

Under the condition
\begin{equation}
\label{allornot3d}
\textrm{rank\ } \mathcal{G}(f_{ab})=\textrm{constant\ }\geq 2 \ (p\in M) \ \Leftrightarrow\ \ 
\textrm{rank\ } \mathcal{G}(\mu_{kl}^j)\stackrel{\eqref{Gramsrel}}{=}\textrm{constant\ }\geq 6 \ (p\in M),
\end{equation}
the (pointwise) extremal solutions form a smooth section of the torsion tensor bundle. If the corresponding linear connection is compatible with the Finsler metric then we have a non-Riemannian generalized Berwald space: recall that condition 
$$\textrm{rank\ } \mathcal{G}(f_{ab})=\textrm{constant\ } \ (p\in M) \ \Leftrightarrow\ \ \textrm{rank\ } \mathcal{G}(\mu_{kl}^j)\stackrel{\eqref{Gramsrel}}{=}\textrm{constant\ } \ (p\in M)$$
obviously holds for a generalized Berwald metric because the parallel transports with respect to a compatible linear connection provide linear isometries between the tangent spaces. If the lower bound in \eqref{allornot3d} does not work then we have a Riemannian space because of Wang's theorem for the maximal dimension of the linear isometry group of a non-Riemannian indicatrix hypersurface (see Corollary \ref{cor:01}):
$$d\leq \frac{(n-1)(n-2)}{2}\ \ \Rightarrow\ \ \textrm{rank\ } \mathcal{G}(f_{ab})=\textrm{constant\ }=\binom{n}{2}-d\geq n-1.$$

Using the orthogonalization process (\ref{ortproc}) we can write that 
$$T(\lambda)=-\sum_{j}\sum_{k < l} \frac{\langle \nu_{kl}^j, h^*\rangle}{\|T(\nu_{kl}^j)\|^2}T(\nu_{kl}^j).$$
Especially, if the rank of the system (\ref{system02}) is maximal then the extremal compatible linear connection is the only compatible linear connection because of
$$\dim W/\textrm{\ Ker}\ T=\dim \textrm{\ Im\ } T=3\binom{3}{2}=\dim V.$$
According to Corollary \ref{cor:01}, it is of zero curvature and its torsion is of the form 
$$T(\lambda)=-\sum_{j=1}^3\sum_{1\leq k < l \leq 3} \frac{\langle \nu_{kl}^j, h^*\rangle}{\|T(\nu_{kl}^j)\|^2}T(\nu_{kl}^j).$$ 

\subsection{The general case} Since any  connection parameter belongs to a triplet of indices, all computations can be repeated in case of spaces of dimension greater than three. Let $j=1, \ldots, n$ be given. We are going to compute the elements of the matrix table in the $j$-th block. Using \eqref{torsionmapping}, it follows that
$$T_{ab}^c(\mu)=\frac{1}{2}\int_{S_p^*} \mu_c f_{ab}+\mu_a f_{cb}-\mu_b f_{ca}\, \mu^*\ \ \Rightarrow\ \ T_{ab}^c(\mu_{kl}^j)=\frac{1}{2}\int_{S_p^*} f_{kl} \left(\delta^j_c f_{ab}+\delta^j_a f_{cb}-\delta^j_b f_{ca}\right)\, \mu^*.$$
Introducing the quantities
$$A_{I(a,b)}=\left\{\begin{array}{ll}
\ \ A_{(a-1)n-\frac{a(a-1)}{2}+b-a}& (a<b)\\
\ \ \ 0 & (a=b)\\
-A_{(b-1)n-\frac{b(b-1)}{2}+a-b}& (a>b)
\end{array}\right.
$$
for the row vectors of the Gram matrix $\mathcal{G}(f_{ab})$, we have vectors of the form
$$\frac{1}{2} \left(\delta^j_c A_{I(a,b)}+\delta^j_a A_{I(c,b)}-\delta^j_b A_{I(c,a)}\right)=\left\{\begin{array}{ll}
\ \ \ 0& (j\neq c, j\neq a, j\neq b, a<b)\\
&\\
\ \ A_{I(a,b)}/2 & (j=c, j\neq a, j\neq b, a<b)\\
&\\
\ \ A_{I(a,j)} & (j=c, j\neq a, j= b, a<j)\\
&\\
\ \ A_{I(j,b)} & (j=c, j= a, j\neq b, j<b)\\
&\\
-A_{I(c,a)}/2 & (j\neq c, j\neq a, j=b, a<j)\\
&\\
\ \ A_{I(c,b)}/2 & (j\neq c, j=a, j\neq b, j<b)
\end{array}\right.$$
in the rows of the $j$-th block of the matrix table. Taking a non-trivial vanishing linear combination of the rows of the matrix table we have the following equation in the $j$-th block:
$$0=\sum_{a=1}^{j-1} \lambda_{aj}^j A_{I(a,j)}+\sum_{b=j+1}^{n}\lambda_{jb}^j A_{I(j,b)}+$$
$$\frac{1}{2}\left(\hat{\sum_{1\leq a < b \leq n}} \lambda_{ab}^j A_{I(a,b)}-\sum_{c\neq j} \sum_{a=1}^{j-1}\lambda_{aj}^c A_{I(c,a)}+\sum_{c\neq j}\sum_{b=j+1}^{n}\lambda_{jb}^c A_{I(c,b)}\right),$$
where the symbol $\hat{\sum}$ means a $j$-free sum. Therefore
$$\textrm{the coefficient of the vector\ }A_{I(\alpha, \beta)}=\left\{\begin{array}{ll}
\lambda_{\alpha, j}^j & (j=\beta)\\
&\\
\lambda_{j\beta}^j& (j=\alpha)\\
&\\
\dfrac{\lambda_{\alpha \beta}^j+\lambda_{j \beta}^{\alpha}-\lambda_{ j \alpha}^{\beta}}{2}& (j < \alpha<\beta)\\
&\\
\dfrac{\lambda_{\alpha \beta}^j+\lambda_{\alpha j}^{ \beta}+\lambda_{ j \beta}^{\alpha}}{2}& (\alpha< j < \beta)\\
&\\
\dfrac{\lambda_{\alpha \beta}^j+\lambda_{\alpha j}^{ \beta}-\lambda_{ \beta j}^{\alpha}}{2}& (\alpha<\beta < j)
\end{array}\right.
$$ 
for any indices $1\leq \alpha < \beta \leq n$. As the block index $j$ runs through the numbers $1, \ldots, n$ we have a corresponding non-trivial vanishing linear combination of the rows of the Gram matrix of system \eqref{system02} (it is a block-diagonal matrix with the Gram matrix $\mathcal{G}(f_{ab})$ along the diagonal $n$ times). By the help of a reasoning similar to the case of dimension three we have the general result.

\section{Extremal compatible linear connections in Randers spaces}

A Finsler metric is called a \emph{Randers metric} if it has the special form
\[ F(x,y) = \alpha(x,y) + \beta(x,y),  \]
where $\alpha(x,y)=\sqrt{\alpha_{ij}(x)y^{i}y^{j}}$ is a norm coming from a Riemannian metric and $\beta$ is a 1-form with the coordinate representation $\beta(x,y)=\beta_j(x)y^j$. Both the metric components of the Riemannian part and the components of the perturbating term are considered on the tangent manifold as composite functions $\alpha_{ij}(x)$ and $\beta_k(x)$, where $x=(x^1, \ldots, x^n)$. It is known \cite{Vin1} that if a linear connection  is compatible to the Randers metric (the parallel transports preserve the Randers length of the tangent vectors), then it is a metric linear connection with respect to the Riemannian part $\alpha$. Therefore $\gamma:=\alpha$ is a convenient choice for a Riemannian environment because $\alpha$ is directly given by the Randers metric. We suppose that $M$ is connected and $\beta_p\neq 0$ at each point of the manifold. Otherwise, if there is a point $p \in M$ of a Randers manifold admitting compatible linear connections such that $\beta_p=0$, then the indicatrix at $p$ is quadratic and the same holds at all points of the (connected) manifold because the indicatrices are related by linear parallel translations. Therefore the metric is Riemannian.

Let a point $p\in M$ be given. In order to make the computations easier, we choose local coordinates around $p$ such that the coordinate vector fields $\partial/\partial u^1, \ldots, \partial/\partial u^n$ form an orthonormal basis at $p\in M$, $\beta_1(p) = \dots = \beta_{n-1}(p)=0$ and $\beta_n(p) \neq 0$, i.e. the coordinate vector fields from $1$ to $n-1$ span the kernel of the linear functional  $\beta_p \colon T_pM \to \mathbb{R}$. Therefore  
\[F(x,y)=\sqrt{\delta_{ij}y^i y^j} +\beta_n(x)y^n = \sqrt{(y^1)^2+\dots+(y^n)^2} + \beta_n(x)y^n \]
and the partial derivatives of $F$ with respect to the vectorial directions are
\[\dfrac{\partial F}{\partial y^k}(x,y)= \dfrac{y^k}{\sqrt{\sum_{i=1}^{n}(y^i)^2}} + \delta^n_k \beta_n(x)\ \ \Rightarrow\ \ y^a \frac{\partial F}{\partial y^b}-y^b\frac{\partial F}{\partial y^a} =\left(y^a \delta^n_b - y^b \delta^n_a \right)\beta_n. \]
Taking $a< b$ we have that 
$$f_{ab}=y^a\delta_b^n \beta_n \ \Rightarrow \ \textrm{rank\ }\mathcal{G}(f_{ab})=n-1$$
because of
$$\int_{S_p^*}y^iy^j\, \mu^*=0\ (i\neq j) \ \ \ \textrm{and}\ \ \ \omega=\int_{S_p^*}(y^1)^2\, \mu^*=\ldots=\int_{S_p^*}(y^n)^2\, \mu^*.$$
The dimension of the linear isometry group of the indicatrix is
$$d=\binom{n}{2}-(n-1)=\binom{n-1}{2}.$$
The linear isometry group is actually the entire orthogonal group of the kernel of the perturbating term. In what follows we are going to prove the characterization theorem \cite{VO} (see also \cite{Vin1}) in the context of the general theory.

\begin{Thm} {\emph{\cite{VO}}} The conditional extremum problem for the torsion of the extremal compatible linear connection of a Randers space is solvable at the point $p\in M$ if and only if 
\begin{equation}
\label{randsolv}
C_{n;i}:=\dfrac{\partial \beta_n}{\partial x^i}-\beta_n \Gamma_{in}^{*n}=0 \ \ \ (i=1, \ldots, n),
\end{equation}
i.e. the dual vector field $\beta^{\sharp}$ of the perturbating term with respect to $\alpha$ is of constant Riemannian length as the point runs through the base manifold. 
\end{Thm}

\begin{Pf} If $h^*=0$ then the uniquely determined solution is the identically zero torsion tensor at the point $p\in M$. Since $\nabla^*$ is the L\'{e}vi-Civita connection of $\alpha$, we have
\begin{equation} \label{randersh}
\begin{gathered}
h_i^*=X_i^{h^*} F=X_i^{h^*}\alpha + X_i^{h^*} \beta = X_i^{h^*} \beta = \frac{\partial \beta_s y^s}{\partial x^i}-y^j {\Gamma}^{*k}_{ij}\circ \pi \frac{\partial \beta_s y^s}{\partial y^k} = \\
\frac{\partial \beta_s}{\partial x^i}y^s-y^j {\Gamma}^{*k}_{ij}\circ \pi  \beta_k = y^j \left(  \frac{\partial \beta_j}{\partial x^i} - \beta_n \Gamma^{*n}_{ij} \circ \pi  \right) =: y^j C_{j;i};
\end{gathered}
\end{equation}
that is, $y^j C_{j;i}=0$. Differentiating with respect to $\partial/\partial y^n$ we have \eqref{randsolv}. Otherwise, by Theorem \ref{solvability03}, the solvability condition is $\mathrm{Ker} \ T\subset {h^*}^{\bot}$. Using that, by \eqref{torsionmapping}, 
\[ \renewcommand\arraystretch{1.8} T_{ab}^c(\mu)=\frac{\beta_n}{2} \cdot \left\{\begin{array}{ll}
\int_{S_p^*}y^a\mu_c+y^c\mu_a\, \mu^* & (a<b=n, c< n)\\
2\int_{S_p^*}y^a\mu_n \, \mu^*& (a<b=n, c= n)\\
\int_{S_p^*}y^a\mu_b-y^b\mu_a \, \mu^*& (a<b<n, c= n)\\
0& (a< b< n, c< n),
\end{array}\right. \]
the elements of $\mathrm{Ker} \ T$ are characterized by
\begin{equation} \label{randersker}
\mu \in \mathrm{Ker} \ T \ \Leftrightarrow \ T_{ab}^c(\mu)= 0 \ \Leftrightarrow \ \int_{S_p^*}y^a\mu_b \, \mu^* =0 \ \  (a<n, b\leq n).
\end{equation}
On the other hand, by \eqref{randersh}, $h_i^*=y^j C_{j;i}$, so the elements of ${h^*}^{\bot}$ are characterized by
\begin{equation}
\label{ortcomp}
\mu \in {h^*}^{\bot} \ \Leftrightarrow \ \sum_{i=1}^n \int_{S_p^*} \mu_i h_i^* \, \mu^*=0  \ \Leftrightarrow \ \sum_{i=1}^n C_{j;i} \int_{S_p^*} y^j \mu_i\, \mu^*=0.
\end{equation}

Let us take an arbitrary $\mu \in \mathrm{Ker} \ T$; by \eqref{ortcomp}, we have $\mu \in {h^*}^{\bot}$ if and only if 
\begin{equation} \label{randerskerort}
\sum_{i=1}^n C_{j;i} \int_{S_p^*} y^j \mu_i\, \mu^* \overset{\eqref{randersker}}{=} \sum_{i=1}^n C_{n;i} \int_{S_p^*} y^n \mu_i\, \mu^*=0.
\end{equation}
\begin{itemize}
\item Supposing that \eqref{randsolv} holds, it immediately implies \eqref{randerskerort}, and thus $\mathrm{Ker} \ T\subset {h^*}^{\bot}$.
\item Supposing $\mathrm{Ker} \ T\subset {h^*}^{\bot}$, we can choose elements
\[ \mu^1=(y^n, 0, \ldots, 0), \ \mu^2=(0, y^n, 0, \ldots, 0), \ \ldots, \ \mu^n=(0, \ldots, 0, y^n) \in \mathrm{Ker} \ T, \]
i.e. $\mu_i^k=\delta_i^ky^n$ (obviously satisfying condition \eqref{randersker}), which, by necessarily satisfying \eqref{randerskerort}, imply for any $k=1,\dots, n$ that
\[ \sum_{i=1}^n C_{n;i} \int_{S_p^*} y^n \mu_i^k\, \mu^*= \sum_{i=1}^n C_{n;i} \int_{S_p^*} y^n \delta_i^ky^n\, \mu^*= C_{n;k} \int_{S_p^*} (y^n)^2 \mu^*= 0 \ \Leftrightarrow \ C_{n;k}=0. \]
\end{itemize}
As a straightforward calculation shows (see \cite{VO}, subsection 4.5) equation \eqref{randsolv} is equivalent to the vanishing of the derivative of the norm square of the dual vector field with respect to the Riemannian metric $\alpha$.
\end{Pf}

\subsection{Computations in 3D} Using that
$$\int_{S_p^*}y^iy^j\, \mu^*=0\ (i\neq j) \ \ \ \textrm{and}\ \ \ \omega=\int_{S_p^*}(y^1)^2\, \mu^*=\ldots=\int_{S_p^*}(y^n)^2\, \mu^*,$$
a straightforward computation shows that the matrix table is of the form

\begin{center}
{\tiny \addtolength{\tabcolsep}{-2pt}{
\begin{tabular}{|c|c|c|c||c|c|c||c|c|c|}
\hline
&&&&&&&&&\\                                                                
$T_{ab}^c$&$T(\mu_{12}^1)$&$T(\mu_{13}^1)$&$T(\mu_{23}^1)$&$T(\mu_{12}^2)$&$T(\mu_{13}^2)$&$T(\mu_{23}^2)$&$T(\mu_{12}^3)$&$T(\mu_{13}^3)$&$T(\mu_{23}^3)$\\
&&&&&&&&&\\
\hline
&&&&&&&&&\\
$T_{12}^1$ &0&0&0&0&0&0&0&0& 0\\
&&&&&&&&&\\
\hline 
&&&&&&&&&\\
$T_{13}^1$ &0&$\omega \beta_3^2$&0&0&0&0&0&0&0\\
&&&&&&&&&\\
\hline
&&&&&&&&&\\
$T_{23}^1$ &0&0&$\omega \beta_3^2/2$&0&$\omega \beta_3^2/2$&0&0&0&0\\
&&&&&&&&&\\
\hline
&&&&&&&&&\\
$T_{12}^2$ &0 &0&0&0&0&0&0&0&0\\
&&&&&&&&&\\
\hline
&&&&&&&&&\\
$T_{13}^2$ &0&0&$\omega \beta_3^2/2$&0&$\omega \beta_3^2/2$&0&0&0&0\\
&&&&&&&&&\\
\hline
&&&&&&&&&\\
$T_{23}^2$ &0 &0&0&0&0&$\omega \beta_3^2$&0&0&0\\
&&&&&&&&&\\
\hline
&&&&&&&&&\\
$T_{12}^3$ &0&0&$-\omega \beta_3^2/2$&0&$\omega \beta_3^2/2$&0&0&0&0\\
&&&&&&&&&\\
\hline
&&&&&&&&&\\
$T_{13}^3$ &0 &0&0&0&0&0&0&$\omega \beta_3^2$&0\\
&&&&&&&&&\\
\hline
&&&&&&&&&\\
$T_{23}^3$ &0 &0&0&0&0&0&0&0&$\omega \beta_3^2$\\
&&&&&&&&&\\
\hline
\end{tabular}}}
\end{center}

and the Gram matrix of the basis $T(\mu_{13}^1),\ T(\mu_{23}^1), \ T(\mu_{13}^2),\ T(\mu_{23}^2), \ T(\mu_{13}^3),\ T(\mu_{23}^3)$ is 
$$
\mathcal{G}(T(\mu_{kl}^j))=\beta_3^4 \omega^2 \begin{bmatrix}
1&0&0&0&0&0\\
0&3/4&1/4&0&0&0\\
0&1/4&3/4&0&0&0\\
0&0&0&1&0&0\\
0&0&0&0&1&0\\
0&0&0&0&0&1
\end{bmatrix}\ (1\leq k < l=3, j=1, 2, 3).
$$
Since the solution of the conditional extremum problem is of the form 
$$T(\lambda)=r_{13}^1 T(\mu_{13}^1)+r_{23}^1 T(\mu_{23}^1)+r_{13}^2 T(\mu_{13}^2)+r_{23}^2 T(\mu_{23}^2)+r_{13}^3 T(\mu_{13}^3)+r_{23}^3 T(\mu_{23}^3),$$
the substitution of the elements $\mu=\mu_{13}^1,\ \mu_{23}^1, \ \mu_{13}^2,\ \mu_{23}^2, \ \mu_{13}^3,\ \mu_{23}^3$ into the constraint equations \eqref{solvability01} gives that
\begin{equation}
\begin{bmatrix}
A_{13}^1\\
A_{23}^1\\
A_{13}^2\\
A_{23}^2\\
A_{13}^3\\
A_{23}^3\\
\end{bmatrix}+\beta_3^4\omega^2\begin{bmatrix}
1&0&0&0&0&0\\
0&3/4&1/4&0&0&0\\
0&1/4&3/4&0&0&0\\
0&0&0&1&0&0\\
0&0&0&0&1&0\\
0&0&0&0&0&1
\end{bmatrix}\begin{bmatrix}
r_{13}^1\\
r_{23}^1\\
r_{13}^2\\
r_{23}^2\\
r_{13}^3\\
r_{23}^3\\
\end{bmatrix}=\begin{bmatrix}
0\\
0\\
0\\
0\\
0\\
0\\
\end{bmatrix},
\end{equation}
where the affine terms are
\begin{equation}
\begin{aligned}
A_{13}^1=\langle h^*, \mu_{13}^1\rangle&=\int_{S_p^*}h^*_1 f_{13}\, \mu^*=\int_{S_p^*}\left(X_1^{h^*}F\right) f_{13}\, \mu^*=\omega\beta_3 \left(\frac{\partial \beta_1}{\partial x^1}-\beta_3\Gamma_{11}^{*3}\right),\\
A_{23}^1=\langle h^*, \mu_{23}^1\rangle&=\int_{S_p^*}h^*_1 f_{23}\, \mu^*=\int_{S_p^*}\left(X_1^{h^*}F\right) f_{23}\, \mu^*=\omega\beta_3 \left(\frac{\partial \beta_2}{\partial x^1}-\beta_3\Gamma_{12}^{*3}\right),\\
A_{13}^2=\langle h^*, \mu_{13}^2\rangle&=\int_{S_p^*}h^*_2 f_{13}\, \mu^*=\int_{S_p^*}\left(X_2^{h^*}F\right) f_{13}\, \mu^*=\omega\beta_3 \left(\frac{\partial \beta_1}{\partial x^2}-\beta_3\Gamma_{21}^{*3}\right),\\
A_{23}^2=\langle h^*, \mu_{23}^2\rangle&=\int_{S_p^*}h^*_2 f_{23}\, \mu^*=\int_{S_p^*}\left(X_2^{h^*}F\right) f_{23}\, \mu^*=\omega\beta_3 \left(\frac{\partial \beta_2}{\partial x^2}-\beta_3\Gamma_{22}^{*3}\right),\\
A_{13}^3=\langle h^*, \mu_{13}^3\rangle&=\int_{S_p^*}h^*_3 f_{13}\, \mu^*=\int_{S_p^*}\left(X_3^{h^*}F\right) f_{13}\, \mu^*=\omega\beta_3 \left(\frac{\partial \beta_1}{\partial x^3}-\beta_3\Gamma_{31}^{*3}\right),\\
A_{23}^3=\langle h^*, \mu_{23}^3\rangle&=\int_{S_p^*}h^*_3 f_{23}\, \mu^*=\int_{S_p^*}\left(X_3^{h^*}F\right) f_{23}\, \mu^*=\omega\beta_3 \left(\frac{\partial \beta_2}{\partial x^3}-\beta_3\Gamma_{32}^{*3}\right).\\
\end{aligned}
\end{equation}
Therefore
\begin{equation}
\begin{aligned}
\begin{bmatrix}
r_{13}^1\\
r_{23}^1\\
r_{13}^2\\
r_{23}^2\\
r_{13}^3\\
r_{23}^3\\
\end{bmatrix}&=-
\frac{1}{\beta_3^4\omega^2}\begin{bmatrix}
1&0&0&0&0&0\\
0&3/4&1/4&0&0&0\\
0&1/4&3/4&0&0&0\\
0&0&0&1&0&0\\
0&0&0&0&1&0\\
0&0&0&0&0&1
\end{bmatrix}^{-1}
\begin{bmatrix}
A_{13}^1\\
A_{23}^1\\
A_{13}^2\\
A_{23}^2\\
A_{13}^3\\
A_{23}^3\\
\end{bmatrix}\\
&=-
\frac{1}{\beta_3^4\omega^2}\begin{bmatrix}
1&\ \ 0&\ \ 0&0&0&0\\
0&\ \ 3/2&-1/2&0&0&0\\
0&-1/2&\ \ 3/2&0&0&0\\
0&\ \ 0&\ \ 0&1&0&0\\
0&\ \ 0&\ \ 0&0&1&0\\
0&\ \ 0&\ \ 0&0&0&1
\end{bmatrix}
\begin{bmatrix}
A_{13}^1\\
A_{23}^1\\
A_{13}^2\\
A_{23}^2\\
A_{13}^3\\
A_{23}^3\\
\end{bmatrix}\\
\end{aligned}
\end{equation}
and
$$T_{12}^1(\lambda)=T_{12}^2(\lambda)=0, \ T_{13}^1(\lambda)=r_{13}^1T_{13}^1(\mu_{13}^1)=-\frac{1}{\beta_3} \left(\frac{\partial \beta_1}{\partial x^1}-\beta_3\Gamma_{11}^{*3}\right),
$$
$$T_{23}^2(\lambda)=r_{23}^2T_{23}^2(\mu_{23}^2)=-\frac{1}{\beta_3} \left(\frac{\partial \beta_2}{\partial x^2}-\beta_3\Gamma_{22}^{*3}\right),$$
$$T_{13}^3(\lambda)=r_{13}^3T_{13}^3(\mu_{13}^3)=-\frac{1}{\beta_3} \left(\frac{\partial \beta_1}{\partial x^3}-\beta_3\Gamma_{13}^{*3}\right),$$
$$T_{23}^3(\lambda)=r_{23}^3T_{23}^3(\mu_{23}^3)=-\frac{1}{\beta_3} \left(\frac{\partial \beta_2}{\partial x^3}-\beta_3\Gamma_{23}^{*3}\right),$$
$$T_{23}^1(\lambda)=r_{23}^1T_{23}^1(\mu_{23}^1)+r_{13}^2T_{23}^1(\mu_{13}^2)=\frac{\beta_3^2 \omega}{2}\left(r_{23}^1+r_{13}^2\right)=\Gamma^{*3}_{12}-\frac{1}{2\beta_3}\left(\frac{\partial \beta_2}{\partial x^1}+\frac{\partial \beta_1}{\partial x^2}\right),$$
$$T_{13}^2(\lambda)=r_{23}^1T_{13}^2(\mu_{23}^1)+r_{13}^2T_{13}^2(\mu_{13}^2)=\frac{\beta_3^2 \omega}{2}\left(r_{23}^1+r_{13}^2\right)=\Gamma^{*3}_{12}-\frac{1}{2\beta_3}\left(\frac{\partial \beta_2}{\partial x^1}+\frac{\partial \beta_1}{\partial x^2}\right),$$
$$T_{12}^3(\lambda)=r_{23}^1T_{12}^3(\mu_{23}^1)+r_{13}^2T_{12}^3(\mu_{13}^2)=-\frac{\beta_3^2 \omega}{2}\left(r_{23}^1-r_{13}^2\right)=\frac{1}{\beta_3}\left(\frac{\partial \beta_2}{\partial x^1}-\frac{\partial \beta_1}{\partial x^2}\right).$$

\section{Summary} The method presented in the paper provides an intrinsic characterization of generalized Berwald spaces. They are beyond Riemannian spaces in the sense that the indicatrices are more general (not necessarily quadratic) hypersurfaces in the tangent spaces. We are close to the Riemannian geometry in the sense that generalized Berwald spaces admit linear connections on the base manifold such that parallel transports preserve the length of the tangent vectors (compatibility condition). The basic problem is the intrinsic characterization of such a linear connection. Our idea is based on the averaged Riemannian metric and the extremal compatible linear connection minimizing the length of its torsion tensor (the bundle metric is given by the averaged Riemannian metric). To find the extremal compatible linear connection of a connected generalized Berwald space we have to solve a conditional extremum problem in the tangent spaces point by point. It is a hybrid problem because the (quadratic) objective function is given on a finite dimensional vector space (the fiber of the torsion tensor bundle at a given point of the base manifold) but the constraint equations are given as functional equations on the indicatrix hypersurface at a given point of the base manifold. First of all we formulate the sufficient and necessary condition of the solvability (Theorem \ref{solvability03}). The Lagrange principle on function spaces allows us to find a general formula (\ref{solution}) for the solution. Beyond the general formula we also give the solution in terms of intrinsic quantities of the space (Theorem \ref{nd}). During the solution of the conditional extremum problem, we present computations in 2D and 3D (Theorem \ref{2d}, Theorem \ref{3d}) to illustrate the steps of the procedure. Finally, we apply the general results to the special class of Randers metrics.

\section{Acknowledgments}

M\'{a}rk Ol\'{a}h has received funding from the HUN-REN Hungarian Research Network.

\end{document}